\documentclass[final,1p,times]{elsarticle}

\usepackage{amssymb}
 \usepackage{amsthm}
\usepackage{amscd}
\usepackage{amsmath}
\usepackage{amsfonts}
\usepackage{amssymb}
\usepackage{graphicx}
\newtheorem{theorem}{Theorem}
\newtheorem{proposition}[theorem]{Proposition}
\newtheorem{lemma}[theorem]{Lemma}
\usepackage{mathrsfs}
\usepackage{titletoc}

\newcommand{\ra}{\rightarrow}
\newcommand{\p}{\partial}
\newcommand{\f}{\frac}

\newcommand{\be}{\begin{equation}}
\renewcommand{\ra}{\rightarrow}
\newcommand{\ee}{\end{equation}}
\newcommand{\bea}{\begin{eqnarray}}
\newcommand{\eea}{\end{eqnarray}}
\newcommand{\bna}{\begin{eqnarray*}}
\newcommand{\ena}{\end{eqnarray*}}

\renewcommand{\le}{\left}
\newcommand{\ri}{\right}

\journal{***}

\begin{document}

\begin{frontmatter}
\title{A Trudinger-Moser inequality on compact Riemannian surface involving Gaussian curvature}

\author{Yunyan Yang}
 \ead{yunyanyang@ruc.edu.cn}
\address{ Department of Mathematics,
Renmin University of China, Beijing 100872, P. R. China}

\begin{abstract}
  Motivated by a recent work of X. Chen and M. Zhu (Commun. Math. Stat., 1 (2013) 369-385),
  we establish a Trudinger-Moser inequality on  compact Riemannian surface without boundary.
  The proof is based on blow-up analysis together with Carleson-Chang's result (Bull. Sci.
  Math. 110 (1986) 113-127).
  This inequality is different from the classical one, which is due to L. Fontana
  (Comment. Math. Helv., 68 (1993) 415-454), since the Gaussian curvature is involved.
  As an application, we improve Chen-Zhu's result as follows: A modified Liouville energy of
  conformal Riemannian metric has a uniform lower bound, provided that the Euler characteristic
  is nonzero and the volume of the conformal surface has a uniform positive lower bound.
\end{abstract}


\end{frontmatter}

\titlecontents{section}[0mm]
                       {\vspace{.2\baselineskip}}
                       {\thecontentslabel~\hspace{.5em}}
                        {}
                        {\dotfill\contentspage[{\makebox[0pt][r]{\thecontentspage}}]}
\titlecontents{subsection}[3mm]
                       {\vspace{.2\baselineskip}}
                       {\thecontentslabel~\hspace{.5em}}
                        {}
                       {\dotfill\contentspage[{\makebox[0pt][r]{\thecontentspage}}]}

\setcounter{tocdepth}{2}


\section{Introduction}

Let $(\Sigma,g)$ be a compact Riemannian surface without boundary, $K_g$ be its Gaussian curvature and $\chi(\Sigma)$ be
its Euler characteristic.
Let $W^{1,2}(\Sigma)$ be the completion of $C^\infty(\Sigma)$ under the norm
\be\label{Sobolev-norm}\|u\|_{W^{1,2}(\Sigma)}=\le(\int_\Sigma(|\nabla_gu|^2+u^2)dv_g\ri)^{1/2},\ee
where $\nabla_g$ is the gradient operator and $dv_g$ is the Riemannian volume element.
As a limit case of the Sobolev embedding theorem, the Trudinger-Moser inequality \cite{24,19,17,22,14} plays an important role in
analysis and geometry. In an elegant paper \cite{Fontana},
L. Fontana  proved that D. Adams' results \cite{Adams},
Trudinger-Moser inequalities for higher order derivatives, still hold on compact Riemannian manifolds without
boundary. Among those, there is the following
\be\label{Fontana-ineq}\sup_{u\in W^{1,2}(\Sigma),\,\int_\Sigma udv_g=0,\,\int_\Sigma|\nabla_gu|^2dv_g\leq 1}
\int_\Sigma e^{\gamma u^2}dv_g<+\infty,\quad\forall \gamma\leq 4\pi.\ee
Motivated by works of Adimurthi-Druet \cite{A-D}, the author \cite{Yang-JFA,Yang-Tran,Yang-MZ} and C. Tintarev \cite{Tint},
we obtained in \cite{Yang-JDE} that for any $\alpha<\lambda_g^\ast(\Sigma)$, there holds
\be\label{Fontana-ineq-improve}\sup_{u\in W^{1,2}(\Sigma),\,\int_\Sigma udv_g=0,\,\int_\Sigma|\nabla_gu|^2dv_g-\alpha
\int_\Sigma u^2dv_g\leq 1}
\int_\Sigma e^{\gamma u^2}dv_g<+\infty,\quad\forall \gamma\leq 4\pi,\ee
and extremal function for this inequality exists. Here $\lambda_g^\ast(\Sigma)$ is the first eigenvalue of the Laplace-Beltrami operator with respect to the mean value  zero condition, namely
\be\label{eig-mean}\lambda_g^\ast(\Sigma)=\inf_{u\in W^{1,2}(\Sigma),\int_\Sigma udv_g=0,\int_\Sigma u^2dv_g=1}\int_\Sigma|\nabla_g u|^2dv_g.\ee
Clearly $\lambda_g^\ast(\Sigma)>0$ and (\ref{Fontana-ineq-improve}) improves (\ref{Fontana-ineq}). As a consequence of
(\ref{Fontana-ineq-improve}), we have a weak form of the Trudinger-Moser inequality. Namely, for any $\alpha<\lambda_g^\ast(\Sigma)$, there exists
some constant $C$ depending only on $(\Sigma,g)$ and $\alpha$ such that for all $u\in W^{1,2}(\Sigma)$, there holds
\be\label{weakform}\int_\Sigma|\nabla_gu|^2dv_g-\alpha\int_\Sigma(u-\overline{u})^2dv_g-16\pi
\ln\int_\Sigma e^udv_g+16\pi \overline{u}\geq -C,\ee
where $\overline{u}=\f{1}{{\rm vol}_g(\Sigma)}\int_\Sigma udv_g$. When $\alpha=0$, the inequality was obtained by Ding-Jost-Li-Wang \cite{DJLW}.

For any conformal metric
$\tilde{g}=e^ug$, where $u\in C^2(\Sigma)$, the Liouville energy of $\tilde{g}$ reads
$$L_g(\tilde{g})=\int_\Sigma\ln\f{\tilde{g}}{g}(R_{\tilde{g}}dv_{\tilde{g}}+R_gdv_g),$$
where $R_g$ and $R_{\tilde{g}}$ are twice the Gaussian curvature $K_g$ and $K_{\tilde{g}}$ respectively.
Since
$$R_{\tilde{g}}=e^{-u}(\Delta_gu+R_g),$$
we have
\be\label{Liouv}L_g(\tilde{g})=\int_\Sigma\le(|\nabla_gu|^2+4K_gu\ri)dv_g.\ee
If $\Sigma$ is a topological two sphere and ${\rm vol}_{\tilde{g}}(\Sigma)={\rm vol}_{g}(\Sigma)=4\pi$,
then it was proved by X. Chen and M. Zhu \cite{Chen-Zhu} that there exists some constant $C$ depending only on $(\Sigma,g)$ such that
\be\label{lowerb}L_g(\tilde{g})\geq -C.\ee
This is a very important issue in the Calabi flow \cite{Chen,Chen-Zhu}.
Note that (\ref{weakform}) can be derived from (\ref{Fontana-ineq-improve}).
One would expect an inequality, which is an analog of (\ref{Fontana-ineq-improve}) and stronger than
(\ref{lowerb}).  To state our results,
we fix several notations. Let us first define a function space
\be\label{space}\mathscr{K}_g=\le\{u\in W^{1,2}(\Sigma): \int_\Sigma K_gudv_g=0\ri\}\ee
and an associate eigenvalue of the Laplace-Beltrami operator

\be\label{eigenvalue}\lambda_g(\Sigma)=\inf_{u\in\mathscr{K}_g,\,\int_\Sigma u^2dv_g=1}\int_\Sigma|\nabla_g u|^2dv_g.\ee
Clearly, $\mathscr{K}_g$ is a closed subspace of $W^{1,2}(\Sigma)$. While unlike $\lambda_g^\ast(\Sigma)$ given as in
(\ref{eig-mean}),
$\lambda_g(\Sigma)$ is not necessarily
nonzero. For example, $\lambda_g(\Sigma)=0$ if $K_g\equiv 0$.
In Lemma \ref{lemma2.1} below, we shall describe a necessary and sufficient condition under which
 $\lambda_g(\Sigma)>0$. If
$\alpha<\lambda_g(\Sigma)$ and $u\in \mathscr{K}_g$, we write
\be\label{norm-alpha}\|u\|_{1,\alpha}=\le(\int_\Sigma|\nabla_gu|^2dv_g-\alpha\int_\Sigma u^2dv_g\ri)^{1/2}.\ee
Clearly $\|\cdot\|_{1,\alpha}$ is an equivalent norm to $\|\cdot\|_{W^{1,2}(\Sigma)}$
defined as in (\ref{Sobolev-norm}) on the function space $\mathscr{K}_g$, provided that
 $\alpha<\lambda_g(\Sigma)$. The first and the most important result in this paper can be stated as follows:

\begin{theorem}\label{Theorem 1}
 Let $(\Sigma,g)$  be a compact Riemannian surface without boundary,
$K_g$ be its Gaussian curvature,
and $\mathscr{K}_g$, $\lambda_g(\Sigma)$ be defined as in (\ref{space}), (\ref{eigenvalue})
respectively. Suppose that  the Euler characteristic $\chi(\Sigma)\not=0$. Then for any $\alpha<\lambda_g(\Sigma)$, there holds
\be\label{Fontana-revised}\sup_{u\in \mathscr{K}_g,\, \|u\|_{1,\alpha}\leq 1}\int_\Sigma e^{4\pi u^2}dv_g<+\infty,\ee
where $\|\cdot\|_{1,\alpha}$ is a norm defined as in (\ref{norm-alpha}).
Moreover, $4\pi$ is the best constant, in other words,
 if $e^{4\pi u^2}$ is replaced by $e^{\gamma u^2}$ in (\ref{Fontana-revised}) for any $\gamma>4\pi$, then
the above supremum is infinity.
\end{theorem}

An extremely interesting case of Theorem \ref{Theorem 1} is $\alpha=0$. If $\lambda_g(\Sigma)>0$, the norm $\|\cdot\|_{1,0}$ is well defined on the
function space $\mathscr{K}_g$, and thus (\ref{Fontana-revised}) is an analog of Fontana's inequality
(\ref{Fontana-ineq}).
Furthermore, the following theorem reveals the relation between the Trudinger-Moser inequality and the
topology of $\Sigma$.

\begin{theorem}\label{Theorem2}
 Let $(\Sigma,g)$  be a compact Riemannian surface without boundary, $\mathscr{K}_g$
 be defined as in (\ref{space}).
Then the Trudinger-Moser inequality
\be\label{T-M-0}\sup_{u\in \mathscr{K}_g,\,
\int_\Sigma|\nabla_gu|^2dv_g\leq 1}\int_\Sigma e^{4\pi u^2}dv_g<+\infty\ee
holds if and only if the Euler characteristic $\chi(\Sigma)\not= 0$.
\end{theorem}

Let us explain the relation between (\ref{T-M-0}) and (\ref{Fontana-ineq}) under the assumption that
$\chi(\Sigma)\not= 0$.
First, we can see that both best constants
of the two inequalities are $4\pi$. Second, the subcritical inequalities in both cases are equivalent. Precisely,
the inequalities
$$\sup_{u\in \mathscr{K}_g,\,
\int_\Sigma|\nabla_gu|^2dv_g\leq 1}\int_\Sigma e^{\gamma u^2}dv_g<+\infty,\quad\forall\gamma<4\pi,$$
holds if and only if (\ref{Fontana-ineq}) holds for all $\gamma<4\pi$. Third, in the critical case,
(\ref{T-M-0})
is independent of (\ref{Fontana-ineq}).

Another interesting problem for the Trudinger-Moser inequality is the existence of extremal functions. Pioneer works in
this direction were due to Carleson-Chang \cite{CC}, M. Struwe \cite{Struwe},
F. Flucher \cite{Flucher}, K. Lin \cite{Lin}, Ding-Jost-Li-Wang \cite{DJLW,DJLW-2}, and Adimurthi-Struwe \cite{Adi-Stru}.
Concerning the extremal functions for (\ref{Fontana-revised}), we have the following:

\begin{theorem}\label{Theorem3}
If the Euler characteristic $\chi(\Sigma)\not=0$, then for any
$\gamma\leq 4\pi$ and $\alpha<\lambda_g(\Sigma)$, where $\lambda_g(\Sigma)$ is defined as in
(\ref{eigenvalue}), the supremum in (\ref{Fontana-revised}) can be attained
by some function $u^\ast\in\mathscr{K}_g$ with $\|u^\ast\|_{1,\alpha}\leq 1$.
\end{theorem}

One would ask what will happen when $\chi(\Sigma)=0$. We talk about this situation briefly.
From Lemma \ref{lemma2.1}
below, we know that $\lambda_g(\Sigma)=0$. By the Young inequality $2ab\leq \epsilon a^2+\epsilon^{-1}b^2$,
$\forall\epsilon>0$, and Fontana's inequality (\ref{Fontana-ineq}), we can prove that
\be\label{1alpha}\sup_{u\in W^{1,2}(\Sigma),\, \|u\|_{1,\alpha}\leq 1}\int_\Sigma e^{4\pi u^2}dv_g<+\infty,\quad\forall\alpha<0.\ee
For details of the proof of (\ref{1alpha}), we refer the reader to \cite{do-Yang,Yang-Zhao}. Thus (\ref{1alpha})
 is weaker than (\ref{Fontana-ineq}). There also holds
\bna \sup_{u\in \mathscr{K}_g,\, \|u\|_{1,\alpha}\leq 1}\int_\Sigma e^{4\pi u^2}dv_g\leq
\sup_{u\in W^{1,2}(\Sigma),\, \|u\|_{1,\alpha}\leq 1}\int_\Sigma e^{4\pi u^2}dv_g.\ena
 Hence, in the case $\chi(\Sigma)=0$, (\ref{Fontana-revised}) is still true and
slightly weaker than (\ref{Fontana-ineq}). It was proved by Y. Li \cite{Lijpde} that the extremal function
for (\ref{1alpha}) exists, while it is open whether or not the extremal function for
(\ref{Fontana-revised}) exists under the assumption $\chi(\Sigma)=0$. This issue will not be discussed here.\\

Finally we apply Theorem \ref{Theorem 1} to Chen-Zhu's problem \cite{Chen-Zhu}. Let $\tilde{g}=
e^ug$ be a metric conformal to $g$, where $u\in C^2(\Sigma)$.
If $\chi(\Sigma)\not=0$, we define a modified Liouville energy of $\tilde{g}$,
 by
\be\label{newfunct}\overline{L}_g(\tilde{g})=\int_\Sigma\le(|\nabla_gu|^2+\f{8}{\chi(\Sigma)}K_gu\ri)dv_g.\ee
In particular, if $\Sigma$ is a topological two sphere, then  $\overline{L}_g(\tilde{g})$ coincides with
the Liouville energy $L_g(\tilde{g})$ defined as in (\ref{Liouv}).
 We denote
\be\label{supremum}C_g(\Sigma)=\sup_{u\in\mathscr{K}_g,\,\int_\Sigma|\nabla_gu|^2dv_g\leq 1}
\int_\Sigma e^{4\pi u^2}dv_g.\ee
The following theorem generalizes Chen-Zhu's result (\ref{lowerb}).

\begin{theorem}\label{Theorem4}
 Let $(\Sigma,g)$ be a compact Riemannian surface without boundary.
Suppose that the Euler characteristic $\chi(\Sigma)\not=0$. For any conformal metric
$\tilde{g}=e^ug$ with $u\in C^2(\Sigma)$, if
${\rm vol}_{\tilde{g}}(\Sigma)\geq \mu{\rm vol}_{g}(\Sigma)$
for some constant $\mu>0$,
then there holds
$$\overline{L}_g(\tilde{g})\geq 16\pi\ln\f{\mu{\rm vol}_{g}(\Sigma)}{C_g(\Sigma)},$$
where $\overline{L}_g(\tilde{g})$ and $C_g(\Sigma)$ are defined as in (\ref{newfunct}) and (\ref{supremum})
respectively.
\end{theorem}

The proof of Theorem \ref{Theorem 1} and Theorem \ref{Theorem3} is based on blow-up analysis. We follow the lines of
\cite{Yang-Tran,Yang-JDE}, and thereby follow closely Y. Li \cite{Lijpde} and Adimurthi-Druet \cite{A-D}.
Earlier works had been done by Ding-Jost-Li-Wang \cite{DJLW,DJLW-2} and Adimuthi-Struwe \cite{Adi-Stru}. Both Theorem \ref{Theorem2}
and Theorem \ref{Theorem4} are consequences of Theorem \ref{Theorem 1}.
The following lemma due to Carleson-Chang \cite{CC} will be used in our analysis.

\begin{lemma}\label{CC-lemma}
    Let $\mathbb{B}$ be the unit disc in $\mathbb{R}^2$.
          Assume $\{v_\epsilon\}_{\epsilon>0}$ is a sequence of functions in $W_0^{1,2}(\mathbb{B})$
          with $\int_{\mathbb{B}}|\nabla v_\epsilon|^2dx=1$. If $|\nabla v_\epsilon|^2dx\rightharpoonup
          \delta_0$ as $\epsilon\ra 0$
          weakly in sense of measure. Then
          $$\limsup_{\epsilon\ra 0}\int_{\mathbb{B}}(e^{4\pi v_\epsilon^2}-1)dx\leq\pi
          e.$$
\end{lemma}

\noindent Another key ingredient in our analysis is the well-known Gauss-Bonnet formula
(see for example \cite{Gallot}, Section 3.J.1, p. 176), namely
\be\label{Gauss-Bonnet}\int_\Sigma K_gdv_g=2\pi\chi(\Sigma).\ee

Throughout this paper, we often denote various constants by the same $C$, also we do not distinguish sequence and subsequence.
The remaining part of this paper is organized as follows:  Theorem \ref{Theorem 1} and  Theorem \ref{Theorem3} are
proved in Section 2; Theorem \ref{Theorem2} and Theorem \ref{Theorem4} are proved in Section 3 and Section 4 respectively.

\section{A Trudinger-Moser inequality involving Gaussian curvature}

In this section, we prove Theorem \ref{Theorem 1} and Theorem \ref{Theorem3}
by using blow-up analysis.  We need several preliminary results before beginning the blow-up procedure.

\begin{lemma}\label{lemma2.1}
$\lambda_g(\Sigma)>0$ if and only if $\chi(\Sigma)\not=0$.
\end{lemma}

\noindent{\it Proof}. Let us first prove that $\lambda_g(\Sigma)$ can be attained.
By definition of $\lambda_g(\Sigma)$ (see (\ref{eigenvalue}) before),
we take $u_j\in\mathscr{K}_g$ such that $\int_\Sigma u_j^2dv_g=1$ and $\int_\Sigma |\nabla_gu_j|^2dv_g\ra \lambda_g(\Sigma)$
as $j\ra+\infty$. Clearly $u_j$ is bounded in $W^{1,2}(\Sigma)$. Thus, up to a subsequence, we can assume
\bna
&&u_j\rightharpoonup u_0\quad{\rm weakly\,\,in}\quad W^{1,2}(\Sigma),\\[1.2ex]
&&u_j\ra u_0\quad{\rm strongly\,\,in}\quad L^2(\Sigma).
\ena
It then follows that
\be\label{u0}\int_\Sigma u_0^2dv_g=1,\quad\int_\Sigma|\nabla_g u_0|^2dv_g\leq \lambda_g(\Sigma).\ee
Since $u_j\in\mathscr{K}_g$, we have
\be\label{uzr}\int_\Sigma K_gu_0dv_g=\lim_{j\ra+\infty}\int_\Sigma K_gu_jdv_g=0.\ee
By (\ref{u0}) and (\ref{uzr}), we have that $u_0\in\mathscr{K}_g$ attains $\lambda_g(\Sigma)$.

If $\chi(\Sigma)\not=0$, we claim that $\lambda_g(\Sigma)>0$. Suppose not. We have $\lambda_g(\Sigma)=0$ and
thus $u_0\equiv C$ for some constant $C$.
In view of (\ref{uzr}), we have by using the Gauss-Bonnet formula (\ref{Gauss-Bonnet}) that
$2\pi C\chi(\Sigma)=0$. Hence $C=0$. This contradicts
$\int_\Sigma u_0^2dv_g=1$. On the contrary, if $\chi(\Sigma)=0$, then the Gauss-Bonnet formula implies that
$u\equiv c\in\mathscr{K}_g$ for any $c\in\mathbb{R}$.
Hence $\lambda_g(\Sigma)=0$. $\hfill\Box$\\

An immediate consequence of Lemma \ref{lemma2.1} is the following:
\begin{lemma}\label{lemma2.2} Suppose that $\chi(\Sigma)\not=0$. Then $\forall u\in\mathscr{K}_g$, we have
$\int_\Sigma u^2dv_g\leq \f{1}{\lambda_g(\Sigma)}\int_\Sigma|\nabla_gu|^2dv_g$.
\end{lemma}

\begin{proposition}\label{Proposition2.2'} Suppose that $\chi(\Sigma)\not=0$. There holds
\be\label{sub-k0}\sup_{u\in\mathscr{K}_g,\,
\int_\Sigma|\nabla_gu|^2dv_g\leq 1}
\int_\Sigma e^{\gamma u^2}dv_g<+\infty,\quad\forall\gamma<4\pi.\ee
\end{proposition}

\noindent{\it Proof}. Take any $u\in\mathscr{K}_g$ with $\int_\Sigma|\nabla_gu|^2dv_g\leq 1$. Denote
$\overline{u}=\f{1}{{\rm vol}_g(\Sigma)}\int_\Sigma udv_g$.
For any fixed $\gamma<4\pi$, we can find some $\gamma_0$, say $\gamma_0=(\gamma+4\pi)/2$, and a constant $C$
depending only on $\gamma$ such that
$$\gamma u^2\leq \gamma_0(u-\overline{u})^2+C\overline{u}^2.$$
By Lemma \ref{lemma2.2}, there exists some constant $C$ depending only on $(\Sigma,g)$ such that
$$\overline{u}^2\leq \f{1}{{\rm vol}_g(\Sigma)}\int_\Sigma u^2dv_g\leq C.$$
Then it follows from Fontana's inequality (\ref{Fontana-ineq}) that
$$\int_\Sigma e^{\gamma u^2}dv_g\leq C\int_\Sigma e^{\gamma_0 (u-\overline{u})^2}dv_g\leq C$$
for some constant $C$ depending only on $(\Sigma,g)$ and $\gamma$. Therefore (\ref{sub-k0}) follows.
 $\hfill\Box$\\

We remark that if $\chi(\Sigma)\not=0$, then Proposition \ref{Proposition2.2'} indicates that Fontana's subcritical inequalities imply (\ref{sub-k0}).
Conversely, assuming (\ref{sub-k0}), then using the same argument as in the above proof we can get (\ref{Fontana-ineq})
for any $\gamma<4\pi$. Therefore, (\ref{sub-k0}) is equivalent to (\ref{Fontana-ineq}) with $\gamma<4\pi$.\\

 In the remaining part of this section, we always assume $\chi(\Sigma)\not=0$ and $\alpha<\lambda_g(\Sigma)$.

\begin{lemma}\label{lemma2.3}
For any $0<\epsilon<4\pi$, there exists some $u_\epsilon\in C^1(\Sigma)\cap \mathscr{K}_g$
with $\|u_\epsilon\|_{1,\alpha}=1$ such that
\be\label{sub-c-1}\int_\Sigma e^{(4\pi-\epsilon)u_\epsilon^2}dv_g=\sup_{u\in\mathscr{K}_g,\,
\|u\|_{1,\alpha}\leq 1}
\int_\Sigma e^{(4\pi-\epsilon)u^2}dv_g.\ee
\end{lemma}

\noindent{\it Proof.}
 For any fixed $0<\epsilon<4\pi$, we choose $u_j\in \mathscr{K}_g$ with $\|u_j\|_{1,\alpha}\leq 1$
 such that as $j\ra+\infty$,
\be\label{sub-2}\int_\Sigma e^{(4\pi-\epsilon)u_j^2}dv_g\ra
\sup_{u\in \mathscr{K}_g,\,
\|u\|_{1,\alpha}\leq 1}
\int_\Sigma e^{(4\pi-\epsilon)u^2}dv_g.\ee
By Lemma \ref{lemma2.2}, $u_j$ is bounded in $W^{1,2}(\Sigma)$.
Then we can assume, up to a subsequence, $u_j\rightharpoonup u_\epsilon$  weakly  in
  $W^{1,2}(\Sigma)$,  $u_j\ra u_\epsilon$  strongly in  $L^2(\Sigma)$, and
  $u_j\ra u_\epsilon$ a.e.  in $\Sigma$.
  As such, we have
  $$\int_\Sigma |\nabla_g u_\epsilon|^2 dv_g\leq \limsup_{j\ra+\infty}\int_\Sigma |\nabla_g u_j|^2 dv_g.$$
  This immediately leads to $\|u_\epsilon\|_{1,\alpha}\leq 1$  and
  $$\int_\Sigma|\nabla_g u_j-\nabla_g u_\epsilon|^2dv_g\leq
  1-\|u_\epsilon\|_{1,\alpha}^2+o_j(1).$$
  Observe
  $$(4\pi-\epsilon)u_j^2\leq (4\pi-\epsilon/2)(u_j-u_\epsilon)^2+32\pi^2\epsilon^{-1}u_\epsilon^2.$$
  It follows
  from the H\"older inequality and Proposition \ref{Proposition2.2'} that
  $e^{(4\pi-\epsilon)u_j^2}$ is bounded in $L^q(\Sigma)$ for some $q>1$. Hence $e^{(4\pi-\epsilon)u_j^2}\ra e^{(4\pi-\epsilon)u_\epsilon^2}$
  strongly in $L^1(\Sigma)$. This together with (\ref{sub-2}) leads to (\ref{sub-c-1}).
  Since $u_j\in\mathscr{K}_g$, we have $u_\epsilon\in\mathscr{K}_g$. It is easy to see that $\|u_\epsilon\|_{1,\alpha}=1$.

  By a straightforward calculation, we know that $u_\epsilon$ satisfies the Euler-Lagrange equation
 \be\label{E-L-1}\le\{
  \begin{array}{lll}
  \Delta_g u_\epsilon-\alpha u_\epsilon=\f{1}{\lambda_\epsilon}u_\epsilon e^{(4\pi-\epsilon)u_\epsilon^2}-
  \mu_\epsilon K_g\\[1.5ex]
  \lambda_\epsilon=\int_\Sigma u_\epsilon^2 e^{(4\pi-\epsilon)u_\epsilon^2}dv_g\\
  [1.5ex] \mu_\epsilon=\f{1}{2\pi\chi(\Sigma)}\le(\f{1}{\lambda_\epsilon}
  \int_\Sigma u_\epsilon e^{(4\pi-\epsilon)u_\epsilon^2}dv_g+\alpha\int_\Sigma u_\epsilon dv_g\ri),
  \end{array}
  \ri.\ee
  where $\Delta_g$ is the Laplace-Beltrami operator. Elliptic estimate implies that $u_\epsilon\in C^1(\Sigma)$. $\hfill\Box$

 \begin{lemma}\label{lemma2.4}
  We have
  \be\label{limit}\lim_{\epsilon\ra 0}\int_\Sigma e^{(4\pi-\epsilon)u_\epsilon^2}dv_g=\sup_{u\in \mathscr{K}_g,\,
\|u\|_{1,\alpha}\leq1}
\int_\Sigma e^{(4\pi-\epsilon)u^2}dv_g.\ee
\end{lemma}

\noindent{\it Proof.} By Lemma \ref{lemma2.3}, $\int_\Sigma e^{(4\pi-\epsilon)u_\epsilon^2}dv_g$ is increasing with respect to
$\epsilon>0$. Hence the limit on the left hand side of (\ref{limit}) does exist, possibly it is infinity.
Noting that for any fixed $u\in\mathscr{K}_g$ with $\|u\|_{1,\alpha}\leq 1$, we have
$$\int_\Sigma e^{4\pi u^2}dv_g=\lim_{\epsilon\ra 0}\int_\Sigma e^{(4\pi-\epsilon)u^2}dv_g\leq \lim_{\epsilon\ra 0}
\int_\Sigma e^{(4\pi-\epsilon)u_\epsilon^2}dv_g,$$
and whence
\be\label{11}\sup_{u\in \mathscr{K}_g,\,\|u\|_{1,\alpha}\leq1}
\int_\Sigma e^{4\pi u^2}dv_g\leq \lim_{\epsilon\ra 0}
\int_\Sigma e^{(4\pi-\epsilon)u_\epsilon^2}dv_g.\ee
On the other hand, it is obvious that
\be\label{22}\lim_{\epsilon\ra 0}
\int_\Sigma e^{(4\pi-\epsilon)u_\epsilon^2}dv_g\leq \sup_{u\in \mathscr{K}_g,\,\|u\|_{1,\alpha}\leq1}
\int_\Sigma e^{4\pi u^2}dv_g.\ee
Combining (\ref{11}) and (\ref{22}), we get (\ref{limit}). $\hfill\Box$\\

We now follow the lines of \cite{Yang-Tran,Yang-JDE}, and thereby follow closely
Y. Li \cite{Lijpde} and Adimurthi-Druet \cite{A-D}. Similar blow-up scheme had been
used by Ding-Jost-Li-Wang \cite{DJLW,DJLW-2} and Adimurthi-Struwe \cite{Adi-Stru}.
  Denote $c_\epsilon=|u_\epsilon(x_\epsilon)|=\max_\Sigma |u_\epsilon|$.  If $c_\epsilon$ is bounded, applying elliptic estimates to
  (\ref{E-L-1}), we already conclude the existence of extremal function. Without loss of generality, we may assume
  $c_\epsilon=u_\epsilon(x_\epsilon)\ra +\infty$ and $x_\epsilon \ra p\in\Sigma$
   as $\epsilon\ra 0$.

  \begin{lemma}\label{concenteation}
    $u_\epsilon\rightharpoonup 0$ weakly in $W^{1,2}(\Sigma)$, $u_\epsilon \ra 0$
strongly in $L^q(\Sigma)$ for all $q\geq 1$, and $|\nabla_g u_\epsilon|^2dv_g\rightharpoonup\delta_{p}$ weakly in sense of measure
as $\epsilon\ra 0$,
where $\delta_{p}$ is the usual Dirac measure centered at $p$.
\end{lemma}

\noindent{\it Proof.}  Since $\|u_\epsilon\|_{1,\alpha}=1$ and $u_\epsilon\in\mathscr{K}_g$, it
follows from Lemma \ref{lemma2.2} that $u_\epsilon$ is bounded in $W^{1,2}(\Sigma)$. Precisely, we have
  \be\label{w12}\int_\Sigma|\nabla_gu_\epsilon|^2dv_g+\int_\Sigma u_\epsilon^2dv_g\leq \f{\lambda_g(\Sigma)+1}
  {\lambda_g(\Sigma)-\alpha}.\ee
 In view of (\ref{w12}), without loss of generality, we can assume $u_\epsilon\rightharpoonup u_0$ weakly in $W^{1,2}(\Sigma)$,
and $u_\epsilon\ra u_0$ strongly in $L^q(\Sigma)$ for all $q\geq 1$. It follows that
\be\label{3}\int_\Sigma|\nabla_g u_\epsilon|^2dv_g=1+\alpha\int_\Sigma u_0^2dv_g+o_\epsilon(1)\ee
and
\be\label{4}\int_\Sigma|\nabla_g (u_\epsilon-u_0)|^2dv_g=1-\int_\Sigma|\nabla_g u_0|^2dv_g+\alpha\int_\Sigma u_0^2dv_g+o_\epsilon(1).\ee
Suppose $u_0\not\equiv0$. In view of (\ref{4}), Proposition \ref{Proposition2.2'} together with the H\"older inequality
implies that $e^{4\pi u_\epsilon^2}$ is bounded in $L^q(\Sigma)$
for any fixed $q$ with $1\leq q<1/(1-\|u_0\|_{1,\alpha}^2)$. Applying
elliptic estimates to (\ref{E-L-1}), we have that $u_\epsilon$ is uniformly bounded in $\Sigma$, which contradicts
$c_\epsilon\ra+\infty$.
Therefore $u_0\equiv0$ and (\ref{3}) becomes
\be\label{5}\int_\Sigma|\nabla_g u_\epsilon|^2dv_g=1+o_\epsilon(1).\ee
Suppose $|\nabla_g u_\epsilon|^2dv_g\rightharpoonup \mu$ in sense of measure. If $\mu\not =\delta_{p}$, then in view of (\ref{5}) and $u_0\equiv0$,
 we can choose sufficiently small $r_0>0$
and a cut-off
function $\phi\in C_0^1({B}_{r_0}(p))$ such that $0\leq\phi\leq 1$, $\phi\equiv 1$ on $\mathbb{B}_{r_0/2}(p)$,
$|\nabla_g\phi_0|\leq 4/r_0$ and
$$\limsup_{\epsilon\ra 0}\int_{B_{r_0}(p)}|\nabla_g(\phi u_\epsilon)|^2dv_g<1.$$
Let $\overline{\phi u}_\epsilon=\f{1}{{\rm vol}_g(\Sigma)}\int_\Sigma \phi u_\epsilon dv_g$.
Since $u_\epsilon\ra 0$ strongly in $L^q(\Sigma)$ for all $q\geq 1$, $\overline{\phi u}_\epsilon\ra 0$ as $\epsilon\ra 0$.
Using Fontana's inequality (\ref{Fontana-ineq}), we conclude that $e^{(4\pi-\epsilon)(\phi u_\epsilon)^2}$ is bounded in
 $L^s(B_{r_0}(p))$ for some $s>1$. Note that $\phi\equiv 1$ on $B_{r_0/2}(p)$.
 Applying elliptic estimates to (\ref{E-L-1}), we have that $u_\epsilon$ is uniformly bounded in ${B}_{r_0/2}(p)$,
 which contradicts $c_\epsilon\ra+\infty$ again.
 Therefore $|\nabla_g u_\epsilon|^2dv_g\rightharpoonup \delta_{p}$, and this completes the proof of the lemma.  $\hfill\Box$

\begin{lemma}\label{lemma2.5}
  $\lambda_\epsilon$ has a positive lower bound and $\mu_\epsilon$ is bounded.
  \end{lemma}

\noindent{\it Proof.} Note that $$\int_\Sigma e^{(4\pi-\epsilon)u_\epsilon^2}dv_g\leq\int_\Sigma
  \le(1+(4\pi-\epsilon)u_\epsilon^2e^{(4\pi-\epsilon)u_\epsilon^2}\ri)dv_g= {\rm vol}_g(\Sigma)+(4\pi-\epsilon)\lambda_\epsilon.$$
  This together with Lemma \ref{lemma2.4} implies that $\lambda_\epsilon$ has a positive lower bound.
  By Lemma \ref{concenteation}, we have that $u_\epsilon$ is bounded in $L^2(\Sigma)$. Moreover, since
  $$\f{1}{\lambda_\epsilon}
  \int_\Sigma |u_\epsilon| e^{(4\pi-\epsilon)u_\epsilon^2}dv_g\leq 1+\f{1}{\lambda_\epsilon}e^{4\pi}{\rm vol}_g(\Sigma),$$
  we conclude that $\mu_\epsilon$ is a bounded sequence. $\hfill\Box$\\

   Let $\exp_{x_\epsilon}$ be the exponential map at $x_\epsilon$ and ${\rm inj}_g(\Sigma)$ be
   the injectivity radius of $(\Sigma,g)$. There exists a $\delta$, $0<\delta<
   {\rm inj}_g(\Sigma)$,
     such that for any $\epsilon>0$, $\exp_{x_\epsilon}$ maps the Euclidean disc $\mathbb{B}_\delta(0)\subset
     \mathbb{R}^2$ centered at the origin with radius $\delta$ onto the geodesic disc $B_{\delta}(x_\epsilon)\subset \Sigma$.
   Let
   \bea\label{scal}
   &&r_\epsilon={\sqrt{\lambda_\epsilon}}{c_\epsilon^{-1}}e^{-(2\pi-\epsilon/2)c_\epsilon^2},\\[1.1ex]\nonumber
   &&\widetilde{g}_\epsilon(x)=(\exp_{x_\epsilon}^\ast g)(r_\epsilon x),\quad \forall x\in \mathbb{B}_{\delta r_\epsilon^{-1}}(0).\eea
   For any fixed $\beta$, $0<\beta<4\pi$, we estimate
   \be\label{r-zero}r_\epsilon^2e^{\beta c_\epsilon^2}=\lambda_\epsilon c_\epsilon^{-2}e^{-(4\pi-\epsilon-\beta)c_\epsilon^2}
   \leq c_\epsilon^{-2}\int_\Sigma u_\epsilon^2 e^{\beta u_\epsilon^2}dv_g\ra 0,\ee
   here we have used Proposition \ref{Proposition2.2'} and Lemma \ref{concenteation}. In particular, $r_\epsilon\ra 0$. This leads to
   \be\label{metric}\widetilde{g}_\epsilon\ra \xi\quad{\rm in}\quad C^2_{\rm loc}(\mathbb{R}^2),\ee
   where $\xi$ denotes the Euclidean metric. Define two sequences of blow-up functions on the Euclidean disc
   $\mathbb{B}_{\delta r_\epsilon^{-1}}(0)$ by
  $$\psi_\epsilon(x)=c_\epsilon^{-1}u_\epsilon(\exp_{x_\epsilon}(r_\epsilon x)),\quad
  \varphi_\epsilon(x)=c_\epsilon({u}_\epsilon(\exp_{x_\epsilon}(r_\epsilon x))-c_\epsilon).$$
  It is first discovered by Adimurthi and M. Struwe \cite{Adi-Stru} that the above function sequences are suitable for this kind of
    problems. By the equation (\ref{E-L-1}), we have on $\mathbb{B}_{\delta r_\epsilon^{-1}}(0)$,
  \bea\label{ps-eqn}
  &&-\Delta_{\widetilde{g}_\epsilon}\psi_\epsilon=\alpha r_\epsilon^2\psi_\epsilon+c_\epsilon^{-2}\psi_\epsilon
  e^{(4\pi-\epsilon)(1+\psi_\epsilon)\varphi_\epsilon}
  -r_\epsilon^2c_\epsilon^{-1}\mu_\epsilon
  \widetilde{K}_g,\\[1.2ex]
  \label{ph-eqn}&&-\Delta_{\widetilde{g}_\epsilon}\varphi_\epsilon=\alpha r_\epsilon^2c_\epsilon^2\psi_\epsilon+
  \psi_\epsilon e^{(4\pi-\epsilon)(1+\psi_\epsilon)\varphi_\epsilon}- r_\epsilon^2c_\epsilon \mu_\epsilon
  \widetilde{K}_g,
  \eea
  where $\widetilde{K}_g(x)=K_g(\exp_{x_\epsilon}(r_\epsilon x))$.
   It is easy to see that $\Delta_{\widetilde{g}_\epsilon}\psi_\epsilon\ra 0$ in $L^\infty_{\rm loc}(\mathbb{R}^2)$, $|\psi_\epsilon|\leq 1$
   and $\psi_\epsilon(0)=1$. Applying elliptic estimates to (\ref{ps-eqn}) and noting (\ref{metric}),
   we have $\psi_\epsilon\ra \psi$ in $C^1_{\rm loc}(\mathbb{R}^2)$, where $\psi$ is a distributional solution to
    $$-\Delta_\xi\psi=0\quad{\rm in}\quad \mathbb{R}^2,\quad|\psi|\leq 1,\quad\psi(0)=1.$$
    Then  the Liouville theorem implies that $\psi\equiv 1$ on $\mathbb{R}^2$.

    In view of (\ref{r-zero}), $\Delta_{\widetilde{g}_\epsilon}\varphi_\epsilon$ is
    bounded in $\mathbb{B}_R$ for any fixed $R>0$. Note also that
   $\varphi_\epsilon(x)\leq 0=\varphi_\epsilon(0)$ for all $x\in\mathbb{B}_{
   \delta r_\epsilon^{-1}}(0)$. Applying elliptic estimates to
   (\ref{ph-eqn}), we have $\varphi_\epsilon\ra \varphi$ in $C^1_{\rm loc}(\mathbb{R}^2)$, where $\varphi$ satisfies
   \be\label{bubble}-\Delta_\xi\varphi=e^{8\pi\varphi}\quad{\rm in}\quad\mathbb{R}^2,\quad
   \varphi(0)=0=\sup_{\mathbb{R}^2}\varphi.
   \ee
   Moreover, we have
   \bea\nonumber
   \int_{\mathbb{B}_R(0)}e^{8\pi\varphi(x)}dx&\leq&\limsup_{\epsilon\ra 0}\int_{\mathbb{B}_R(0)}e^{(4\pi-\epsilon)
   (u^2_\epsilon(\exp_{x_\epsilon}(r_\epsilon x))-c_\epsilon^2)}dx\\\nonumber
   &=&\limsup_{\epsilon\ra 0}\int_{\mathbb{B}_{Rr_\epsilon}(0)}e^{(4\pi-\epsilon)
   (u^2_\epsilon(\exp_{x_\epsilon}(y))-c_\epsilon^2)}r_\epsilon^{-2}dy\\\nonumber
   &=&\limsup_{\epsilon\ra 0}\f{1}{\lambda_\epsilon}\int_{B_{Rr_\epsilon}(x_\epsilon)}u_\epsilon^2
   e^{(4\pi-\epsilon)u_\epsilon^2}dv_g\\
   &\leq&1\label{int}
   \eea
   by using (\ref{scal}), change of variables, and $\psi_\epsilon\ra 1$ in $C^1_{\rm loc}(\mathbb{R}^2)$.
   In view of (\ref{bubble}) and (\ref{int}), a result of Chen and Li \cite{CL} implies that
   \be\label{phi} \varphi(x)=-\f{1}{4\pi}\log(1+\pi|x|^2),\quad \forall x\in \mathbb{R}^2.\ee
   As a consequence
   \be\label{bu-1}\int_{\mathbb{R}^2}e^{8\pi\varphi}dx=1.\ee
   In conclusion, we obtain the following:

   \begin{proposition}\label{converg}
   $\psi_\epsilon\ra 1$ in $C^1_{\rm loc}(\mathbb{R}^2)$ and $\varphi_\epsilon\ra \varphi$ in $C^1_{\rm loc}(\mathbb{R}^2)$,
   where $\varphi$ satisfies (\ref{phi}) and (\ref{bu-1}).
   \end{proposition}

   Proposition \ref{converg} provides the convergence behavior of $u_\epsilon$ near the blow-up point $p$.
   For the convergence behavior of $u_\epsilon$ away from $p$, we have the following:

  \begin{proposition}\label{outer}
  $c_\epsilon u_\epsilon\rightharpoonup G$ weakly in $W^{1,q}(\Sigma)$ for all $1<q<2$, and
   $c_\epsilon u_\epsilon\ra G$ in $C^1_{\rm loc}(\Sigma\setminus\{p\})\cap L^2(\Sigma)$,
   where $G$ is a Green function satisfying
   \be\label{Green-1-1}
     \left\{\begin{array}{lll}
            \Delta_g G-\alpha G=\delta_{p}-\frac{1+\alpha\int_\Sigma Gdv_g}{2\pi \chi(\Sigma)}K_g\quad\rm{in}\quad
            \Sigma\\[1.5ex]
            \int_\Sigma G K_gdv_g=0.
         \end{array}\right.\ee
  Moreover, $G$ can be decomposed as
  \be\label{gr-3}G=-\f{1}{2\pi}f(r)\log r+A_p+\psi_\alpha,\ee
  where $r$ denotes the geodesic distance from $p$, $f(r)$ is a nonnegative smooth decreasing function,
   which is equal to $1$ in $B_{{\rm inj}_g(\Sigma)/2}(p)$, and to zero for $r\geq {\rm inj}_g(\Sigma)$,
   $A_p$ is a constant real number, $\psi_\alpha\in C^1(\Sigma)$ with
  $\psi_\alpha(p)=0$.
  \end{proposition}

   \noindent{\it Proof.} With a slight modification of proofs of
   (\cite{Yang-Tran}, Lemmas 4.5-4.9), we obtain
   $c_\epsilon u_\epsilon\rightharpoonup G$ weakly in $W^{1,q}(\Sigma)$ for all $1<q<2$, and
   $c_\epsilon u_\epsilon\ra G$ in $C^1_{\rm loc}(\Sigma\setminus\{p\})\cap L^2(\Sigma)$,
   where $G$ is a distributional solution to (\ref{Green-1-1}).
   It is known (\cite{Aubin}, Section 4.10, p. 106) that there exists some function $h\in L^\infty(\Sigma)$ such that
   $$\Delta_g\le(-\f{1}{2\pi}f(r)\log r\ri)=\delta_p+h$$
  in the distributional sense. Hence
  \be\label{distr}\Delta_g\le(G+\f{1}{2\pi}f(r)\log r\ri)=\alpha G-h-\frac{1+\alpha\int_\Sigma Gdv_g}{2\pi \chi(\Sigma)}K_g
  \ee
  in the distributional sense. Since $G\in L^s(\Sigma)$ for any $s\geq 1$ by the Sobolev embedding theorems,
  the terms on the right hand side of (\ref{distr}) belong to $L^s(\Sigma)$ for all $s\geq 1$.
   By elliptic estimates, $G+\f{1}{2\pi}f(r)\log r\in C^1(\Sigma)$, which implies (\ref{gr-3}). $\hfill\Box$\\

  In the following, we shall derive an upper bound of the integrals $\int_\Sigma e^{(4\pi-\epsilon)u_\epsilon^2}dv_g$. There are two ways to
  obtain the upper bound:
  One is to use the capacity estimate which is due to Y. Li \cite{Lijpde};
  The other is to employ Carleson-Chang's estimate (Lemma \ref{CC-lemma}), which was first used by Li-Liu-Yang \cite{Li-Liu-Yang}.
  Here we prefer to the second way. In view of (\ref{gr-3}), we have
          \bna
          \int_{\Sigma\setminus B_\delta(p)}|\nabla_g
          G|^2dv_g&=&\alpha\int_{\Sigma\setminus
          B_\delta(p)}G^2dv_g-\int_{\p B_\delta(p)}G\f{\p G}{\p \nu}ds_g\\
          &&-\f{1+\alpha\int_\Sigma Gdv_g}{2\pi \chi(\Sigma)}\int_{\Sigma\setminus B_\delta(p)}
          K_g Gdv_g\\
          &=&\f{1}{2\pi}\log\f{1}{\delta}+A_{p}+\alpha\|G\|_2^2+o_\delta(1).
          \ena
          Hence we obtain
          \be\label{Omega-B}
          \int_{\Sigma\setminus B_\delta(p)}|\nabla_g
          u_\epsilon|^2dv_g=\f{1}{c_\epsilon^2}
          \le(\f{1}{2\pi}\log\f{1}{\delta}+A_{p}+\alpha\|G\|_2^2+o_\delta(1)+o_\epsilon(1)\ri).
          \ee
          Let $s_\epsilon=\sup_{\p B_\delta(p)}u_\epsilon$ and
          $\widetilde{u}_\epsilon=(u_\epsilon-s_\epsilon)^+$. Then
          $\widetilde{u}_\epsilon\in W_0^{1,2}(B_\delta(p))$.
          By (\ref{Omega-B}) and the fact that
          $$\int_{B_\delta(p)}|\nabla_g
          u_\epsilon|^2dv_g=1-\int_{\Sigma\setminus B_\delta(p)}|\nabla_g
          u_\epsilon|^2dv_g+\alpha\int_\Sigma u_\epsilon^2dv_g,$$ we have
          $$
          \int_{B_\delta(p)}|\nabla_g
          \widetilde{u}_\epsilon|^2dv_g\leq \tau_\epsilon=1-\f{1}{c_\epsilon^2}
          \le(\f{1}{2\pi}\log\f{1}{\delta}+A_{p}+o_\delta(1)+o_\epsilon(1)\ri).
          $$
          Now we choose an isothermal coordinate system $(U,\phi;\{x^1,x^2\})$ near $p$ such that
          $B_{2\delta}(p)\subset U$,
          $\phi(p)=0$, and the metric $g=e^h(d{x^1}^2+d{x^2}^2)$ for some function $h\in C^1(\phi(U))$
          with $h(0)=0$. Clearly, for any $\delta>0$, there exists some $c(\delta)>0$ with $c(\delta)\ra 0$ as
          $\delta\ra 0$ such that $dv_g\leq (1+c(\delta))dx$ and $\phi(B_\delta(p))\subset
          \mathbb{B}_{\delta(1+c(\delta))}(0)\subset\mathbb{R}^2$.
          Noting that $\widetilde{u}_\epsilon=0$ outside $B_\delta(p)$, we have
          $$\int_{\mathbb{B}_{\delta(1+c(\delta))}(0)}
          |\nabla(\widetilde{u}_\epsilon\circ \phi^{-1})|^2dx=\int_{\phi^{-1}(\mathbb{B}_{\delta(1+c(\delta))}(0))}
          |\nabla_g\widetilde{u}_\epsilon|^2dv_g=\int_{B_\delta(p)}|\nabla_g\widetilde{u}_\epsilon|dv_g\leq \tau_\epsilon.
          $$
          This together with Lemma \ref{CC-lemma}  leads to
          \bea\nonumber
          \limsup_{\epsilon\ra 0}\int_{B_{\delta}(p)}(e^{4\pi
          \widetilde{u}_\epsilon^2/\tau_\epsilon}-1)dv_g&=&\limsup_{\epsilon\ra 0}(1+c(\delta))\int_{\mathbb{B}_{\delta(1+c(\delta))}(0)}(e^{4\pi
          (\widetilde{u}_\epsilon\circ\phi^{-1})^2/\tau_\epsilon}-1)dx\\\label{B-delta}
          &\leq&
          \pi\delta^2(1+c(\delta))^3e.
          \eea
          Note that $|u_\epsilon|\leq c_\epsilon$ and $u_\epsilon/c_\epsilon=1+o_\epsilon(1)$ on the geodesic ball
          $B_{Rr_\epsilon}(x_\epsilon)\subset{\Sigma}$. We estimate on $B_{Rr_\epsilon}(x_\epsilon)$,
          \bna
          (4\pi-\epsilon)
          u_\epsilon^2&\leq&4\pi(\widetilde{u}_\epsilon+s_\epsilon)^2\\
          &\leq&
          4\pi\widetilde{u}_\epsilon^2+8\pi
          s_\epsilon\widetilde{u}_\epsilon+o_\epsilon(1)\\
          &\leq&4\pi\widetilde{u}_\epsilon^2-4\log\delta+8\pi
          A_{p}+o_\epsilon(1)+o_\delta(1)\\
          &\leq&4\pi\widetilde{u}_\epsilon^2/\tau_\epsilon-2\log\delta+4\pi
          A_{p}+o(1).
          \ena
          Therefore
          \bea\nonumber
          \int_{B_{Rr_\epsilon}(x_\epsilon)}e^{(4\pi-\epsilon)
          u_\epsilon^2}dv_g&\leq& \delta^{-2}e^{4\pi
          A_{p}+o(1)}\int_{B_{Rr_\epsilon}(x_\epsilon)}e^{4\pi\widetilde{u}_\epsilon^2/\tau_\epsilon}dv_g\\
          \nonumber
          &=&\delta^{-2}e^{4\pi
          A_{p}+o(1)}\int_{B_{Rr_\epsilon}(x_\epsilon)}(e^{4\pi\widetilde{u}_\epsilon^2/\tau_\epsilon}-1)dv_g+o(1)\\
          \label{ii}
          &\leq&\delta^{-2}e^{4\pi
          A_{p}+o(1)}\int_{B_\delta(p)}(e^{4\pi\widetilde{u}_\epsilon^2/\tau_\epsilon}-1)dv_g+o(1),
          \eea
          where $o(1)\ra 0$ as $\epsilon\ra 0$ first and then $\delta \ra 0$.
          Combining (\ref{B-delta}) with (\ref{ii}), letting $\epsilon\ra 0$ first, and then letting $\delta\ra 0$,
          we conclude
          \be\label{u-p}
          \limsup_{\epsilon\ra 0}\int_{B_{Rr_\epsilon}(x_\epsilon)}e^{(4\pi-\epsilon)
          u_\epsilon^2}dv_g\leq \pi e^{1+4\pi A_{p}}.
          \ee
      By a change of variables and (\ref{bu-1}), there holds
   \bna
   \int_{B_{Rr_\epsilon}(x_\epsilon)}e^{(4\pi-\epsilon)u_\epsilon^2}dv_g&=&(1+o_\epsilon(1))
   \int_{\mathbb{B}_{Rr_\epsilon}(0)}e^{(4\pi-\epsilon)({u}_\epsilon\circ\exp_{x_\epsilon})^{2}}dx\\
   [1.2ex]&=&(1+o_\epsilon(1))
   \int_{\mathbb{B}_{R}(0)}e^{(4\pi-\epsilon)({u}_\epsilon\circ\exp_{x_\epsilon}(r_\epsilon x))^2}r_\epsilon^2dx\\
   &=&(1+o_\epsilon(1))\f{\lambda_\epsilon}{c_\epsilon^2}\int_{\mathbb{B}_{R}(0)}e^{8\pi\varphi}dx\\
   &=&(1+o(1))\f{\lambda_\epsilon}{c_\epsilon^2},
   \ena
   where $o(1)\ra 0$ as $\epsilon\ra 0$ first, and then $R\ra+\infty$.
   This together with (\ref{bu-1}) implies
   \be\label{I1}\lim_{R\ra+\infty}\limsup_{\epsilon\ra 0}\int_{B_{Rr_\epsilon}(x_\epsilon)} e^{(4\pi-\epsilon) u_\epsilon^2}dv_g=
   \limsup_{\epsilon\ra 0}\f{\lambda_\epsilon}{c_\epsilon^2}.\ee
   Using the same argument as (\cite{Lijpde}, Lemma 3.5) (see also \cite{Yang-Tran}, Lemma 3.6), we have
   \be\label{II}\lim_{\epsilon\ra 0}\int_\Sigma e^{(4\pi-\epsilon)u_\epsilon^2}dv_g={\rm vol}_g(\Sigma)+
   \limsup_{\epsilon\ra 0}\f{\lambda_\epsilon}{c_\epsilon^2}.\ee
   Combining (\ref{u-p}), (\ref{I1}) and (\ref{II}), we conclude the following:

   \begin{proposition}\label{upperbound}
   Under the assumption of $c_\epsilon=\max_\Sigma |u_\epsilon|\ra+\infty$, there holds
   $$\sup_{u\in \mathscr{K}_g,\,\|u\|_{1,\alpha}\leq 1}
\int_\Sigma e^{4\pi u^2}dv_g=\limsup_{\epsilon\ra 0}\int_\Sigma e^{(4\pi-\epsilon) u_\epsilon^2}dv_g
\leq {\rm vol}_g(\Sigma)+\pi e^{1+4\pi A_p}.$$
   \end{proposition}

   {\it Completion of the proof of Theorem \ref{Theorem 1}.} If $c_\epsilon$ is bounded, then
   the integral $\int_\Sigma e^{(4\pi-\epsilon)u_\epsilon}dv_g$ is also bounded. Hence Lemma \ref{lemma2.4} leads to
   \be\label{8}\sup_{u\in \mathscr{K}_g,\,\|u\|_{1,\alpha}\leq 1}
\int_\Sigma e^{4\pi u^2}dv_g<+\infty.\ee
While if $c_\epsilon\ra+\infty$, (\ref{8}) follows from Proposition \ref{upperbound} immediately.

Now we prove that $4\pi$ is the best constant. Fixing a point $p\in\Sigma$, we let
$r=r(x)={\rm dist}_g(x,p)$ be the distance from $p$ to $x$ and $B_s=B_s(p)$ be the geodesic ball centered at $p$ with radius $s$.
For sufficiently small $\epsilon>0$, we define a sequence of functions
$$M_\epsilon(x)=\le\{
\begin{array}{lll}
\sqrt{\f{1}{2\pi}\log\f{1}{\epsilon}}, &x\in B_{\epsilon^2}\\[1.5ex]
\f{1}{\sqrt{2\pi\log\f{1}{\epsilon}}}\log\f{\epsilon}{r}, &x\in B_{\epsilon}\setminus B_{\epsilon^2}
\\[1.5ex]
0,&x\in \Sigma\setminus B_{\epsilon}.
\end{array}\ri.$$
Clearly $\|\nabla_gM_\epsilon\|_2=1+O(\epsilon)$. Let $M_\epsilon^\ast=M_\epsilon/\|\nabla_gM_\epsilon\|_2$.
Then we have $\|\nabla_gM_\epsilon^\ast\|_2=1$ and
\be\label{O}
\widetilde{(M_\epsilon^\ast)}_g= \f{\int_\Sigma K_gM_\epsilon^\ast dv_g}{\int_\Sigma K_g dv_g}
=O\le(\epsilon\sqrt{-\log\epsilon}\ri).
\ee
Hence $M_\epsilon^\ast-\widetilde{(M_\epsilon^\ast)}_g\in \mathscr{K}_g$ and
$\|M_\epsilon^\ast-\widetilde{(M_\epsilon^\ast)}_g\|_{1,\alpha}=1+O(\epsilon^2\sqrt{\log\f{1}{\epsilon}})$. If $\gamma>4\pi$, then we take some $\nu$
such that $(1-\nu)\f{\gamma}{2\pi}>2$. We have by (\ref{O}) and an inequality $2ab\leq \nu a^2+b^2/\nu$,
\bna
\int_\Sigma e^{\gamma (M_\epsilon^\ast-\widetilde{(M_\epsilon^\ast)}_g)^2\|M_\epsilon^\ast-\widetilde{(M_\epsilon^\ast)}_g\|_{1,\alpha}^{-2}}
dv_g&\geq&\int_\Sigma e^{\gamma {M_\epsilon^\ast}^2-2\gamma M_\epsilon^\ast\widetilde{(M_\epsilon^\ast)}_g+o_\epsilon(1)}dv_g\\
&\geq&
\int_{B_{\epsilon^2}}e^{(1-\nu)\f{\gamma}{2\pi}\log\f{1}{\epsilon}-\f{\gamma}{\nu}\widetilde{(M_\epsilon^\ast)}_g^2+o_\epsilon(1)}dv_g\\
&=&(1+o_\epsilon(1))\pi\epsilon^{2-(1-\nu)\f{\gamma}{2\pi}}\\
&\ra&+\infty\quad{\rm as}\quad \epsilon\ra 0.
\ena
Therefore for any $\gamma>4\pi$ and $\alpha<\lambda_g(\Sigma)$, there holds
\bna
\sup_{u\in \mathscr{K}_g,\, \|u\|_{1,\alpha}\leq 1}\int_\Sigma e^{\gamma u^2}dv_g\geq
\sup_{\epsilon>0}\int_\Sigma e^{\gamma (M_\epsilon^\ast-\widetilde{(M_\epsilon^\ast)}_g)^2\|M_\epsilon^\ast-\widetilde{(M_\epsilon^\ast)}_g\|_{1,\alpha}^{-2}}
dv_g=+\infty.\ena
This completes the proof of Theorem \ref{Theorem 1}. $\hfill\Box$\\

{\it Completion of the proof of Theorem \ref{Theorem3}.} We shall construct a function sequence
  $\phi_\epsilon$ satisfying
   \be\label{nrm}\int_\Sigma|\nabla_g\phi_\epsilon|^2dv_g-\alpha\int_\Sigma (\phi_\epsilon-\widetilde{({\phi}_\epsilon)}_g)^2dv_g=1\ee
   and
   \be\label{gr}\int_\Sigma e^{4\pi (\phi_\epsilon-\widetilde{({\phi}_\epsilon)}_g)^2}dv_g>{\rm vol}(\Sigma)+\pi e^{1+4\pi A_p}\ee
   for sufficiently small $\epsilon>0$, where
   $$\widetilde{({\phi}_\epsilon)}_g=\f{1}{\int_\Sigma K_gdv_g}\int_\Sigma K_g\phi_\epsilon dv_g.$$
   If there exists such a sequence $\phi_\epsilon$, then we have by Proposition \ref{upperbound}
    that $c_\epsilon$ must be bounded. Applying elliptic estimates
   to (\ref{E-L-1}), we conclude the existence of the desired extremal function.

   Now we construct $\phi_\epsilon$ verifying (\ref{nrm}) and (\ref{gr}).
   Set
     $$\phi_\epsilon=\le\{
     \begin{array}{llll}
     &c+\f{-\f{1}{4\pi}\log(1+\pi\f{r^2}{\epsilon^2})+B}{c}
     \quad &{\rm for} &r\leq R\epsilon\\[1.5ex]
     &\f{G-\eta \psi_\alpha}{c}\quad &{\rm for} & R\epsilon<
     r<2R\epsilon\\[1.2ex]
     &\f{G}{c}\quad &{\rm for} & r\geq 2R\epsilon
     \end{array},
     \ri.$$
     where $\psi_\alpha$ is given by (\ref{gr-3}), $r$ denotes the geodesic distance from $p$, $R=-\log\epsilon$, $\eta\in C_0^\infty(B_{2R\epsilon}(p))$ verifying that $\eta=1$ on $B_{R\epsilon}(p)$ and
     $\|\nabla_g \eta\|_{L^\infty}
     =O(\f{1}{R\epsilon})$, $B$ and $c$ are two constants depending only on $\epsilon$ to be determined
     later.
     In order to assure that $\phi_\epsilon\in W^{1,2}(\Sigma),$ we set
     $$
     c+\f{1}{c}\le(-\f{1}{4\pi}\log(1+\pi R^2)+B\ri)
     =\f{1}{c}\le(-\f{1}{2\pi}\log (R\epsilon)+A_{p}\ri),
     $$
     which gives
     \be\label{2pic2}
     2\pi c^2=-\log\epsilon-2\pi B+2\pi A_{p}+\f{1}{2}\log \pi
     +O(\f{1}{R^2}).
     \ee
     Note that $\int_\Sigma K_gGdv_g=0$. We estimate
     \bea\nonumber
     \int_{\Sigma\setminus B_{R\epsilon}(p)}|\nabla_g G|^2dv_g&=&-\int_{\Sigma\setminus B_{R\epsilon}(p)}G\Delta_gGdv_g+\int_{\p (\Sigma\setminus B_{R\epsilon}(p))}G\f{\p G}{\p\nu}ds\\ \nonumber
     &=&\alpha\int_{\Sigma\setminus B_{R\epsilon}(p)}G^2dv_g-\f{1+\alpha\int_\Sigma Gdv_g}{2\pi\chi(\Sigma)}\int_{\Sigma\setminus B_{R\epsilon}(p)}K_gGdv_g\\
     \nonumber &&-\int_{\p B_{R\epsilon}(p)}G\f{\p G}{\p\nu}ds\\\label{nabl-G}
     &=&-\f{1}{2\pi}\log(R\epsilon)+\alpha\|G\|_2^2+A_p+O(R\epsilon\log(R\epsilon)).
     \eea
     Since $\psi_\alpha\in C^1(\Sigma)$ and $\psi_\alpha(p)=0$, we have
     \bea\label{nabl-eta}&&\int_{B_{2R\epsilon}\setminus B_{R\epsilon}(p)}|\nabla_g\eta|^2\psi_\alpha^2dv_g=O((R\epsilon)^2),\\
     \label{Geta}&&\int_{B_{2R\epsilon}\setminus B_{R\epsilon}(p)}\nabla_g G\nabla_g\eta\psi_\alpha dv_g=O((R\epsilon)^2),\\
     \label{rep}&&\int_{B_{R\epsilon}(p)}|\nabla_g\phi_\epsilon|^2dv_g=\f{1}{c^2}\le(\f{1}{2\pi}\log R+\f{\log\pi}
     {4\pi}-\f{1}{4\pi}+O(\f{1}{R^2})\ri).\eea
     Combining (\ref{nabl-G})-(\ref{rep}), we obtain
     \bea\nonumber
     \int_{\Sigma}|\nabla_g \phi_\epsilon|^2dv_g&=&\f{1}{4\pi c^2}\le(
     2\log\f{1}{\epsilon}+\log\pi-1+4\pi A_{p}+4\pi\alpha\|G\|_2^2\ri.\\\label{gad}
     &&\le.+O(\f{1}{R^2})+
     O(R\epsilon\log(R\epsilon))\ri).
     \eea
     Observing
     \bna
     \int_\Sigma K_g\phi_\epsilon dv_g&=&\f{1}{c}\le(\int_{\Sigma\setminus B_{2R\epsilon}(p)}K_gGdv_g+O(R\epsilon\log(R\epsilon))\ri)\\
     &=&\f{1}{c}\le(-\int_{B_{2R\epsilon}(p)}K_gGdv_g+O(R\epsilon\log(R\epsilon))\ri)\\
     &=&\f{1}{c}O(R\epsilon\log(R\epsilon)),
     \ena
     we have $\widetilde{({\phi}_\epsilon)}_g=\f{1}{c}O(R\epsilon\log(R\epsilon))$. Hence
     \bna\nonumber
     \int_\Sigma(\phi_\epsilon-\widetilde{({\phi}_\epsilon)}_g)^2dv_g&=&\int_\Sigma\phi_\epsilon^2dv_g+\widetilde{({\phi}_\epsilon)}_g^2
     {\rm vol}_g(\Sigma)-2\widetilde{({\phi}_\epsilon)}_g\int_\Sigma\phi_\epsilon dv_g\\
     \label{l2}&=&\f{1}{c^2}\le(\int_\Sigma G^2dv_g+O(R\epsilon\log(R\epsilon))\ri).
     \ena
     This together with (\ref{gad}) yields
     \bna \|\phi_\epsilon-\widetilde{({\phi}_\epsilon)}_g\|_{1,\alpha}^2&=&\int_\Sigma|\nabla_g\phi_\epsilon|^2dv_g-\alpha
     \int_\Sigma(\phi_\epsilon-\widetilde{({\phi}_\epsilon)}_g)^2dv_g\\
     &=&
     \f{1}{4\pi c^2}\le(
     2\log\f{1}{\epsilon}+\log\pi-1+4\pi A_{p}
     +O(\f{1}{R^2})+
     O(R\epsilon\log(R\epsilon))\ri).\ena
     Let $\phi_\epsilon$ satisfy (\ref{nrm}), i.e. $\|\phi_\epsilon-\widetilde{({\phi}_\epsilon)}_g\|_{1,\alpha}=1$. Then we have
     \be\label{c2}
     c^2=-\f{\log\epsilon}{2\pi}+\f{\log\pi}{4\pi}-\f{1}{4\pi}+A_{p}
     +O(\f{1}{R^2})+O(R\epsilon\log(R\epsilon)).
     \ee
     It follows from (\ref{2pic2}) and (\ref{c2}) that
     \be{\label{B}}
     B=\f{1}{4\pi}+O(\f{1}{R^2})+O(R\epsilon\log(R\epsilon)).
     \ee
     Clearly we have on $B_{R\epsilon}(p)$
     $$4\pi(\phi_\epsilon-\widetilde{({\phi}_\epsilon)}_g)^2\geq 4\pi
     c^2-2\log(1+\pi\f{r^2}{\epsilon^2})+8\pi B+O(R\epsilon\log(R\epsilon)).$$
     This together with (\ref{c2}) and (\ref{B}) yields
     \be\label{BRE}
     \int_{B_{R\epsilon}(p)} e^{4\pi(\phi_\epsilon-\widetilde{({\phi}_\epsilon)}_g)^2}dv_g\geq
     \pi e^{1+4\pi A_{p}}
     +O(\f{1}{(\log\epsilon)^2}).
     \ee
     On the other hand,
     \bea
      \label{O-BRE}\int_{\Sigma\setminus
      B_{R\epsilon}(p)}e^{4\pi(\phi_\epsilon-\widetilde{({\phi}_\epsilon)}_g)^2}dv_g&\geq&\int_{\Sigma\setminus
      B_{2R\epsilon}(p)}\le(1+4\pi(\phi_\epsilon-\widetilde{({\phi}_\epsilon)}_g)^2\ri)dv_g{\nonumber}\\
      &\geq& {\rm vol}_g(\Sigma)+4\pi\f{\|G\|_2^2}{c^2}+o(\f{1}{c^2}).
     \eea
     Recalling (\ref{c2}) and combining (\ref{BRE}) and (\ref{O-BRE}),
     we conclude (\ref{gr}) for sufficiently small $\epsilon>0$. This completes the proof of Theorem \ref{Theorem3}.
     $\hfill\Box$

     \section{The Trudinger-Moser inequality versus the topology}

     In this section, we prove Theorem \ref{Theorem2}. The proof is based on Theorem \ref{Theorem 1} and the
     Gauss-Bonnet formula (\ref{Gauss-Bonnet}).\\

     {\it Proof of Theorem \ref{Theorem2}.}
     If $\chi(\Sigma)\not=0$, then we have $\lambda_g(\Sigma)>0$ by Lemma \ref{lemma2.1}. Hence Theorem \ref{Theorem 1}
     holds, in particular, (\ref{Fontana-revised}) holds for $\alpha=0$. This is exactly (\ref{T-M-0}).

     Now we show that if (\ref{T-M-0}) holds, then $\chi(\Sigma)$ must be nonzero. Suppose on the contrary $\chi(\Sigma)=0$.
     It follows from the Gauss-Bonnet formula (\ref{Gauss-Bonnet}) that
     $$\int_\Sigma K_gdv_g=2\pi\chi(\Sigma)=0.$$
     Hence every constant function $u_k=k$ $(k\in\mathbb{N})$ satisfies
     $\int_\Sigma K_gu_kdv_g=0$
     and
     $\int_\Sigma|\nabla_gu_k|^2dv_g=0$.
     Thus $u_k\in \mathscr{K}_g$ and $\|\nabla_gu_k\|_2\leq 1$. Whence, by (\ref{T-M-0}), we have that
     \be\label{bn}\sup_k\int_\Sigma e^{4\pi u_k^2}dv_g\leq C\ee
     for some constant $C$ depending only on $(\Sigma,g)$. But we also have
     $$\int_\Sigma e^{4\pi u_k^2}dv_g=\int_\Sigma e^{4\pi k^2}dv_g\ra +\infty,\quad{\rm as}\quad
     k\ra+\infty.$$
     This contradicts (\ref{bn}). Therefore $\chi(\Sigma)\not=0$, and the proof of Theorem \ref{Theorem2} is finished. $\hfill\Box$

     \section{Lower bound of the modified Liouville energy}

     In this section, as in the proof of Theorem \ref{Theorem2}, we prove Theorem \ref{Theorem4} by combining
     Theorem \ref{Theorem 1} and the
     Gauss-Bonnet formula (\ref{Gauss-Bonnet}).\\

     {\it Proof of Theorem \ref{Theorem4}.}
     Denote
     $$C_{g,\alpha}(\Sigma)=\sup_{u\in\mathscr{K}_g,\,\|u\|_{1,\alpha}\leq 1}\int_\Sigma e^{4\pi u^2}dv_g.$$
     Since $\chi(\Sigma)\not=0$, by
     Lemma \ref{lemma2.1} and Theorem \ref{Theorem 1}, we have $\lambda_g(\Sigma)>0$ and
     $C_{g,\alpha}(\Sigma)<+\infty$ for any $\alpha<\lambda_g(\Sigma)$.
     Let $u\in W^{1,2}(\Sigma)$ be such that $\int_\Sigma|\nabla_gu|^2dv_g\not=0$.
     In view of the Gauss-Bonnet formula (\ref{Gauss-Bonnet}), we set
     $$\widetilde{u}_g=\f{\int_\Sigma K_gudv_g}{\int_\Sigma K_gdv_g}=\f{1}{2\pi\chi(\Sigma)}\int_\Sigma K_gudv_g.$$
     Clearly $u-\widetilde{u}_g\in\mathscr{K}_g$ and $\|u-\widetilde{u}_g\|_{1,\alpha}>0$ for any
     $\alpha<\lambda_g(\Sigma)$.   By the Young
     inequality $2ab\leq \epsilon a^2+\epsilon^{-1}b^2$, we have
     \bna
     u&=&\f{u-\widetilde{u}_g}{\|u-\widetilde{u}_g\|_{1,\alpha}}\|u-\widetilde{u}_g\|_{1,\alpha}+\widetilde{u}_g\\
     &\leq& 4\pi\f{(u-\widetilde{u}_g)^2}{\|u-\widetilde{u}_g\|_{1,\alpha}^2}+
     \f{1}{16\pi}\|u-\widetilde{u}_g\|_{1,\alpha}^2+\widetilde{u}_g.
     \ena
     This leads to
     $$\ln \int_\Sigma e^udv_g\leq \ln\int_\Sigma e^{4\pi\f{(u-\widetilde{u}_g)^2}{\|u-\widetilde{u}_g\|_{1,\alpha}^2}}dv_g
     +\f{1}{16\pi}\|u-\widetilde{u}_g\|_{1,\alpha}^2+\widetilde{u}_g.$$
     Using the assumption $\tilde{g}=e^ug$ and ${\rm vol}_{\tilde{g}}(\Sigma)=\int_\Sigma e^udv_g\geq \mu{\rm vol}_g(\Sigma)$,
     we obtain
     \bna\|u-\widetilde{u}_g\|_{1,\alpha}^2+16\pi \widetilde{u}_g&\geq& -16\pi\ln C_{g,\alpha}(\Sigma)+16\pi\ln
     \int_\Sigma e^udv_g\\&\geq& -16\pi\ln C_{g,\alpha}(\Sigma)+16\pi\ln(\mu {\rm vol_g(\Sigma)}),\ena
     or equivalently
     $$\int_\Sigma|\nabla_gu|^2dv_g-\alpha\int_\Sigma|u-\widetilde{u}_g|^2dv_g+\f{8}{\chi(\Sigma)}\int_\Sigma K_gudv_g
     \geq  16\pi\ln\f{\mu{\rm vol}_{g}(\Sigma)}{C_{g,\alpha}(\Sigma)}.$$
     In particular, choosing $\alpha=0$ in the above inequality, we have
     $$\overline{L}_g(\tilde{g})=\int_\Sigma|\nabla_gu|^2dv_g+\f{8}{\chi(\Sigma)}\int_\Sigma K_gudv_g
     \geq 16\pi\ln\f{\mu{\rm vol}_{g}(\Sigma)}{C_{g,0}(\Sigma)}.$$
     Noting that $C_{g,0}(\Sigma)=C_{g}(\Sigma)$ defined as in (\ref{supremum}), we finish the proof of Theorem \ref{Theorem4}.
     $\hfill\Box$\\

 {\bf Acknowledgement.} This work is supported by the National Science Foundation of China (Grant No.11171347 and Grant
 No. 11471014).

\end{document}